\documentclass[11pt,a4paper,leqno]{amsart}
\usepackage[latin1]{inputenc}
\usepackage[T1]{fontenc}
\usepackage{amsfonts}
\usepackage{amsmath}
\usepackage{amssymb}
\usepackage{eurosym}
\usepackage{mathrsfs}
\usepackage{palatino}
\usepackage{color}
\usepackage{esint}
\usepackage{url}
\usepackage{verbatim}
\usepackage{graphicx}
\usepackage{enumerate}
\usepackage[pagebackref,hypertexnames=false, colorlinks, citecolor=blue, linkcolor=blue, urlcolor=red]{hyperref}

\allowdisplaybreaks[4]

\numberwithin{equation}{section}

\swapnumbers
\theoremstyle{plain}
\newtheorem{thm}[equation]{Theorem}

\newtheorem{prop}[equation]{Proposition}
\newtheorem{cor}[equation]{Corollary}
\theoremstyle{definition}
\newtheorem{defn}[equation]{Definition}

\theoremstyle{remark}
\newtheorem{rem}[equation]{Remark}

 \email{zynn@tju.edu.cn}
 \email{fandashan2@zjnu.edu.cn}
\makeatletter
\@namedef{subjclassname@2020}{  \textup{2020} Mathematics Subject Classification}
\makeatother
\newcommand{\func}[1]{\operatorname{#1}}

\begin{document}
\title{A multilinear oscillatory integral along moment curves}
\author{Yunan Zeng}
\author{Dashan Fan}
\address[Y.Z.]{Center for Applied Mathematics, Tianjin University, Weijin
Road 92, 300072 Tianjin, China}
\address[D.F.]{Department of Mathematical Science, University of
Wisconsin-Milwaukee, Milwaukee, WI 53201, USA}

\begin{abstract}
For moment curves \(\gamma_d\), we study the truncations
\[
T_{\beta,\gamma_d,\varepsilon}(\vec f)(x)=\int_{\varepsilon<|t|<1}e^{i|t|^{-\beta}}
\prod_{i=1}^{d}f_i\bigl(x-\gamma_{d,i}(t)\bigr)\,\frac{dt}{|t|},\qquad \beta>0
\]
We establish variational estimates for this operator.
The proof uses PET induction to obtain power decay on dyadic annuli.
\end{abstract}

\maketitle

\allowdisplaybreaks


\section{Introduction}

Let $d\geq 2$ and let $0<m_{1}<\cdots <m_{d}$ be rational numbers. Fix signs 
$\sigma _{1},\ldots ,\sigma _{d}\in \{-1,1\}$ and define $\gamma
_{d}(t)=(\gamma _{d,1}(t),\ldots ,\gamma _{d,d}(t))$ by 
\begin{equation*}
\gamma _{d,i}(t)=%
\begin{cases}
t^{m_{i}}, & t\geq 0, \\ 
\sigma _{i}|t|^{m_{i}}, & t<0.%
\end{cases}%
\end{equation*}%
We define the multilinear oscillatory integral along moment curve $\gamma
_{d}(t)$ by 
\begin{equation*}
T_{\beta ,\gamma _{d}}(f_{1},\ldots ,f_{d}):=\int_{-1}^{1}e^{i|t|^{-\beta
}}\prod_{i=1}^{d}f_{i}(x-\gamma _{d,i}(t))\frac{dt}{|t|},\qquad \func{where}%
\,\beta >0.
\end{equation*}

Multilinear singular integrals along curves form a central object in modern
harmonic analysis. The model case is the bilinear Hilbert transform, whose $%
L^{p}$ boundedness was obtained by Lacey--Thiele \cite{LT}, and which has driven an extensive theory of
curve-supported and surface-supported multilinear operators. The companion
oscillatory problem, where the kernel is multiplied by a highly oscillatory
phase such as $e^{i\left\vert t\right\vert ^{-\beta }},$ was initiated by
Carbery--Christ--Wright \cite{CCW}, who showed that the oscillation $e^{i\left\vert t\right\vert
^{-\beta }}$ improves the mapping properties of singular integral families,
a phenomenon later strengthened and extended in numerous works. For the
multilinear setting along moment curves, the recent breakthrough of
Guo--Li--Zhan \cite{GLZ} established the uniform single-scale boundedness of
multilinear Hilbert transforms along moment curves (see Theorem \ref{thm GLZ}). We
observe that Theorem \ref{thm GLZ} removes the classical "non-degeneracy" restrictions
on the exponents and provides the uniform, polynomial-independent estimates
that make a variational treatment of the oscillatory problem feasible.
Motivated by these observations, we establish a q-variation estimate for the
family of truncations of $T_{\beta ,\gamma _{d}},$ which is strictly
stronger than mere boundedness: it simultaneously yields the $L^{r}$
boundedness of the principal value operator, quantitative control of the
maximal truncation, and a priori information on the pointwise convergence of 
$T_{\beta ,\gamma _{d},\varepsilon }$ \ as $\varepsilon \rightarrow 0$. Our
first main result is the following theorem.

\begin{thm}
\label{main thm} Let $1<p_{i}<\infty $ for $1\leq i\leq d$, $r>1/2$ and $%
\sum_{i=1}^{d}1/p_{i}=1/r$. Then for every $\beta >0$ and $q>1$, we have 
\begin{equation*}
\Big\|V_{q}\big(T_{\beta ,\gamma _{d},\varepsilon }(f_{1},\ldots
,f_{d}):0<\varepsilon <1\big)\Big\|_{L^{r}(\mathbb{R})}\lesssim
\prod_{i=1}^{d}\Vert f_{i}\Vert _{L^{p_{i}}(\mathbb{R}).}
\end{equation*}%
Here, for $0<\varepsilon <1$, the truncated operator is given by 
\begin{equation*}
T_{\beta ,\gamma _{d},\varepsilon }(\vec{f})(x):=\int_{\varepsilon
<|t|<1}e^{i|t|^{-\beta }}\prod_{i=1}^{d}f_{i}\bigl(x-\gamma _{d,i}(t)\bigr)\,%
\frac{dt}{|t|},
\end{equation*}%
and, for $q\geq 1,$ the $q$-variation of the family of truncations is
defined by 
\begin{equation*}
V_{q}\bigl(T_{\beta ,\gamma _{d},\varepsilon }(\vec{f}):0<\varepsilon <1%
\bigr)(x):=\sup_{0<\varepsilon _{0}<\cdots <\varepsilon _{N}<1}\Big(%
\sum_{j=1}^{N}\bigl|T_{\beta ,\gamma _{d},\varepsilon _{j}}(\vec{f}%
)(x)-T_{\beta ,\gamma _{d},\varepsilon _{j-1}}(\vec{f})(x)\bigr|^{q}\Big)%
^{1/q}.
\end{equation*}
\end{thm}

As an immediate consequence of Theorem \ref{main thm} and the basic properties of $%
V_{q} $ we have:

\begin{cor}
Let $1<p_i<\infty$ for $1\leq i\leq d$, $r>1/2$, and $%
\sum_{i=1}^{d}1/p_i=1/r $. Then, for every $\beta>0$, the following
estimates hold.

\begin{enumerate}
\item \textit{Boundedness of the operator.} 
\begin{equation*}
\bigl\|T_{\beta ,\gamma _{d}}(f_{1},\ldots ,f_{d})\bigr\|_{L^{r}(\mathbb{R}%
)}\lesssim \prod_{i=1}^{d}\Vert f_{i}\Vert _{L^{p_{i}}(\mathbb{R})}.
\end{equation*}

\item \textit{Boundedness of the maximal truncation.} 
\begin{equation*}
\bigl\|T_{\beta ,\gamma _{d}}^{\ast }(f_{1},\ldots ,f_{d})\bigr\|_{L^{r}(%
\mathbb{R})}\lesssim \prod_{i=1}^{d}\Vert f_{i}\Vert _{L^{p_{i}}(\mathbb{R}%
)},
\end{equation*}%
where 
\begin{equation*}
T_{\beta ,\gamma _{d}}^{\ast }(f_{1},\ldots ,f_{d})(x):=\sup_{0<\varepsilon
<1}\bigl|T_{\beta ,\gamma _{d},\varepsilon }(f_{1},\ldots ,f_{d})(x)\bigr|.
\end{equation*}
\end{enumerate}
\end{cor}

A key feature of the proof is that it does not rely on the exact monomial
structure of $\gamma _{d}$ beyond a single non-vanishing order of contact at
the origin. With a slight modification, the argument extends to arbitrary
real polynomial curves:

\begin{cor}
\label{cor:polynomial-curves} Let $P_1,\ldots,P_d$ be real polynomials. Let $%
1<p_i<\infty$ and $1/r=\sum_{i=1}^d1/p_i$, with $r>1$. For every $\beta>0$
and $q>1$, we have 
\begin{equation*}
\bigl\|V_q(T_{\beta,\vec P,\varepsilon}(\vec f):0<\varepsilon<1)\bigr\|_{L^r(%
\mathbb{R})}\lesssim \prod_{i=1}^d\|f_i\|_{L^{p_i}(\mathbb{R})},
\end{equation*}
where 
\begin{equation*}
T_{\beta,\vec P,\varepsilon}(\vec
f)(x)=\int_{\varepsilon<|t|<1}e^{i|t|^{-\beta}}\prod_{i=1}^d f_i(x-P_i(t))\,%
\frac{dt}{|t|}.
\end{equation*}
\end{cor}

\begin{rem}
If an analogue of Theorem~\ref{thm GLZ} held for the polynomial shifts
uniformly over the dyadic annuli, then the proof of Corollary~\ref%
{cor:polynomial-curves} would also cover $1/2<r\leq 1$. Such a uniform
estimate can fail below $r=1$, as the example $d=2$, $P_{1}(t)=t$, $%
P_{2}(t)=t+t^{2}$ shows; it remains open to determine conditions on the
polynomials $P_{1},\ldots ,P_{d}$ under which the variational estimate holds
for $1/2<r\leq 1$.
\end{rem}

The proof of the main theorem relies on three ingredients. (i) Oscillatory
decay. By adapting the PET induction of Bergelson-Leibman to our oscillatory
weight and applying van der Corput, we obtain a power saving $\lambda
^{-c_{0}},$ that is uniform across all relevant scales and parameters
(Proposition \ref{prop PET}). (ii) Single-scale decay. Interpolating this decay with
the uniform boundedness of Guo-Li-Zhan yields a single-scale maximal
estimate carrying a decaying factor $\lambda ^{-\eta }$ (Proposition \ref{prop Sb}).
(iii) Variational summation. We upgrade this pointwise decay to a $V_{q}$
bound, sum the geometrically decaying contributions over dyadic annuli, and
extend the result to rational exponents (Theorem \ref{main thm}).

In Section 2 we develop the calculus of multiplicative differences,
the power-decay Proposition \ref{prop PET} (PET induction and van der Corput),
the uniform single-scale bound (Theorem \ref{thm GLZ}),
and the single-scale maximal estimate with decay (Proposition \ref{prop Sb}). Section
3 proves Theorem \ref{main thm} and Corollary \ref{cor:polynomial-curves}: the decay of Proposition \ref{prop Sb} is
turned into a variational estimate, summed over dyadic annuli to conclude
the integer case, extended to rational exponents, and finally generalized to
polynomial curves.

\bigskip

\section{A Local Oscillatory Estimate}

\bigskip

In this section we establish the single-scale maximal estimate with decay,
Proposition 2.4, which is the main input for the variational argument in
Section 3. Throughout, we fix $d\geq 2$, rational exponents $%
0<m_{1}<m_{2}<\cdot \cdot \cdot <m_{d},$ signs $\sigma _{1},\ldots ,\sigma
_{d}\in \{-1,1\}$ and $\beta >0;$ implicit constants may depend on these and
on $p_{i},r,q,$ but are independent of the scale, the annulus location, and $%
\varepsilon .$

\begin{defn}
Let $w:\mathbb{R}\to\mathbb{C}$ be a measurable function. For $h\in\mathbb{R}
$, define the multiplicative difference of $w$ by 
\begin{equation*}
w_h(u):=w(u+h)\overline{w(u)}.
\end{equation*}
For $K\geq1$ and $\vec h=(h_1,\ldots,h_K)\in\mathbb{R}^K$, define $w_{\vec
h} $ recursively by 
\begin{equation*}
w_{h_1,\ldots,h_k}(u):=w_{h_1,\ldots,h_{k-1}}(u+h_k)\,\overline{%
w_{h_1,\ldots,h_{k-1}}(u)},\qquad 1\leq k\leq K.
\end{equation*}
\end{defn}

The proof of the following proposition uses PET induction, with the
differencing step in \cite[Proposition~4.41]{wan} adapted to the oscillatory
weight $w$. In the study of multiple polynomial ergodic averages, PET
induction is a standard tool for reducing the complexity of polynomial
families. See \cite{BergelsonPET,BL,CFH} for classical applications and \cite%
{KMPW,KMT} for more recent developments.

\begin{prop}
\label{prop PET} Let $\beta>0$ and let $1\leq m_1<\cdots<m_d$ be integers.
For $j\in\mathbb{Z}_{\geq0}$ and any choice of signs $\tau_1,\ldots,\tau_d%
\in\{-1,1\}$, set 
\begin{equation*}
P_i(u)=\tau_i\,2^{-(m_i-m_1)j}u^{m_i},\qquad 1\leq i\leq d.
\end{equation*}
For $k\in\mathbb{Z}$, let $Q=[k,k+1]$ and $Q^*=[k-1,k+2]$. Then, for every $%
\lambda\geq1$, every subinterval $J\subset[1/2,1]$, and every $g_i\in
L^\infty_{\mathrm{loc}}(\mathbb{R})$, we have 
\begin{equation*}
\bigg\|\int_J e^{i\lambda u^{-\beta}}\prod_{i=1}^d g_i\bigl(\,\cdot-P_i(u)%
\bigr)\,\frac{du}{u} \bigg\|_{L^1(Q)}\leq
C\lambda^{-c_0}\prod_{i=1}^d\|g_i\|_{L^\infty(Q^*)}.
\end{equation*}
Here $c_0>0$ depends only on $d$ and $m_d$, and $C>0$ depends only on $\beta$
and $(m_1,\ldots,m_d)$. In particular, the estimate is uniform in $\lambda$, 
$j$, $\tau$, $k$, and $J$.
\end{prop}

\begin{proof}
We first normalize the inputs. Without loss of generality, we may assume
that $\Vert g_{i}\Vert _{L^{\infty }(Q^{\ast })}\neq 0$. Define 
\begin{equation*}
G_{i}(y)=\frac{g_{i}(k+y)}{\Vert g_{i}\Vert _{L^{\infty }(Q^{\ast })}}%
1_{[-1,2]}(y)\qquad \func{for}1\leq i\leq d.
\end{equation*}%
By a change of variables, we only need to prove 
\begin{equation*}
\bigg\|\int_{J}e^{i\lambda u^{-\beta }}\prod_{i=1}^{d}G_{i}\bigl(\,\cdot
-P_{i}(u)\bigr)\,\frac{du}{u}\bigg\|_{L^{1}([0,1])}\leq C\lambda ^{-c_{0}}.
\end{equation*}%
Duality gives 
\begin{equation*}
\begin{aligned} &\bigg\|\int_J e^{i\lambda u^{-\beta}}\prod_{i=1}^d
G_i\bigl(\,\cdot-P_i(u)\bigr)\,\frac{du}{u}\bigg\|_{L^1([0,1])}\\
&\qquad\qquad=\sup_{\|G_0\|_{L^\infty([0,1])}\leq
1}\Big|\int_{0}^{1}G_0(x)\int_{J}e^{i\lambda
u^{-\beta}}\prod_{i=1}^{d}G_i(x-P_i(u))\frac{du}{u}dx\Big| \end{aligned}
\end{equation*}%
We extend $G_{0}$ by zero from $[0,1]$ to $[-1,3)$ and $G_{1},\ldots ,G_{d}$
by zero from $[-1,2]$ to $[-1,3)$, and then extend all these functions $4$%
-periodically to $\mathbb{R}$. Let $\mathbb{T}=\mathbb{R}/4\mathbb{Z}$ and
let $\mu $ be its normalized Haar measure, given by $d\mu (x)=dx/4$ on the
fundamental interval $[-1,3)$. Set 
\begin{equation*}
\Lambda :=\int_{\mathbb{T}}G_{0}(x)\int_{J}e^{i\lambda u^{-\beta
}}\prod_{i=1}^{d}G_{i}(x-P_{i}(u))\frac{du}{u}d\mu (x).
\end{equation*}%
Thus, it suffices to prove that 
\begin{equation*}
|\Lambda |\lesssim \lambda ^{-c_{0}}.
\end{equation*}%
We now apply PET differencing. The Cauchy-Schwarz inequality gives that 
\begin{equation*}
\begin{aligned} |\Lambda|^2&\leq \int_{\mathbb{T}}\Big|\int_J e^{i\lambda
u^{-\beta}}\prod_{i=1}^{d}G_i(x-P_i(u))\frac{du}{u}\Big|^2d\mu(x)\\
&\leq\int_{-\frac{1}{2}}^{\frac{1}{2}}
\Big|\int_{\mathbb{T}}\int_{\mathbb{R}}w(u+h)\overline{w(u)}%
\prod_{i=1}^{d}G_i(x-P_i(u+h))\overline{G_i(x-P_i(u))}dud\mu(x)\Big|dh
\end{aligned}
\end{equation*}%
where $w(u)=\frac{1_{J}(u)}{u}e^{i\lambda u^{-\beta }}$.

Choose a nonconstant polynomial $P_j$ of minimal degree among $%
P_1,\ldots,P_d $. In the spatial integral, substitute $x+P_j(u)$ for $x$ and
let $G^{\prime }_0(x)$ be the product of all factors independent of $u$ with 
$|G^\prime_0(x)|\leq1$. Then 
\begin{equation*}
|\Lambda|^2\leq\int_{-1/2}^{1/2}\Bigg|\int_{\mathbb{T}} G^{\prime
}_0(x)\int_{\mathbb{R}} w_h(u)\prod_{\ell=1}^{d^{\prime }} G^{\prime }_\ell%
\bigl(x-P^{\prime }_\ell(u)\bigr)\,du\,d\mu(x)\Bigg|\,dh,
\end{equation*}
where $d^{\prime }\leq 2d-1$ and $(P^{\prime }_\ell)_{\ell=1}^{d^{\prime }}$
consists of all nonconstant terms in the two families 
\begin{equation*}
P_i(u)-P_j(u)\qquad\text{and}\qquad P_i(u+h)-P_j(u),\qquad 1\leq i\leq d,
\end{equation*}
with each occurrence retained. The corresponding function $G^{\prime }_\ell$
is $\overline{G_i}$ for a term from the first family and $G_i$ for a term
from the second. Since $G^{\prime }_0$ is the product of all factors
independent of $u$, we have $\|G^{\prime }_\ell\|_{L^\infty(\mathbb{T})}\leq
1$ for every $0\leq \ell\leq d^{\prime }$.

The resulting integral has the same form as the original one, so the
differencing step can be repeated. After at most $K$ rounds, all polynomial
shifts become constant. Applying Cauchy--Schwarz in the shift variables at
each step, we obtain that 
\begin{equation*}
|\Lambda |^{2^{K}}\leq \int_{\lbrack -1/2,1/2]^{K}}\Big|\int_{\mathbb{R}}w_{%
\vec{h}}(u)\,du\Big|\,d\vec{h},
\end{equation*}%
where $\vec{h}=(h_{1},\ldots ,h_{K})$ and $K$ depends only on $d$, $m_{d}$.
Applying the van der Corput estimate, we have 
\begin{equation*}
\Big|\int_{\mathbb{R}}w_{\vec{h}}(u)\,du\Big|\lesssim _{\beta ,K}\min \{1,%
\frac{1}{\lambda \prod_{l=1}^{K}|h_{l}|}\}\leq \lambda
^{-1/2}\prod_{l=1}^{K}|h_{l}|^{-1/2}.
\end{equation*}%
Combining the preceding estimates, we obtain $|\Lambda |\lesssim \lambda
^{-c_{0}}$. Here $c_{0}=2^{-K-1}$.
\end{proof}

Before proving the next proposition, we recall a uniform single-scale
estimate of Guo, Li and Zhan.

\begin{thm}[\protect\cite{GLZ}, Theorem~3.2]
\label{thm GLZ} Let $d\geq2$ and let $m_1,\ldots,m_d$ be distinct positive
integers. Let $\rho\in C_c^\infty(\mathbb{R})$ be supported away from the
origin, and let $\alpha_1,\ldots,\alpha_d$ be nonzero real numbers. Suppose
that $1\leq p_i<\infty$ for $1\leq i\leq d, r>1/2$ and $1/r=\sum_{i=1}^d1/{%
p_i}.$ Then 
\begin{equation*}
\bigg\|\int_{\mathbb{R}}\rho(u)\prod_{i=1}^d g_i(\,\cdot-\alpha_i
u^{m_i})\,du\bigg\|_{L^r(\mathbb{R})}\lesssim\prod_{i=1}^d\|g_i\|_{L^{p_i}(%
\mathbb{R})}.
\end{equation*}
The implicit constant is independent of $\alpha_1,\ldots,\alpha_d$.
\end{thm}

Combining Proposition~\ref{prop PET} with Theorem~\ref{thm GLZ}, we obtain
the following single-scale maximal estimate with decay.

\begin{prop}
\label{prop Sb} For $b\in \lbrack 1/2,1]$, define 
\begin{equation*}
S_{b}(g_{1},\ldots ,g_{d})(x):=\int_{1/2}^{b}e^{i\lambda u^{-\beta
}}\prod_{i=1}^{d}g_{i}(x-P_{i}(u))\,\frac{du}{u},
\end{equation*}%
where the polynomials $P_{i}$ are as in Proposition~\ref{prop PET}. Let $%
1<p_{i}<\infty $ for $1\leq i\leq d$, $1/r=\sum_{i=1}^{d}1/p_{i}$ and $r>1/2$%
. Then there exist $\eta >0$ and $C>0$ such that 
\begin{equation*}
\Big\|\sup_{1/2\leq b\leq 1}\big|S_{b}(g_{1},\ldots ,g_{d})\big|\Big\|%
_{L^{r}(\mathbb{R})}\leq C\lambda ^{-\eta }\prod_{i=1}^{d}\Vert g_{i}\Vert
_{L^{p_{i}}(\mathbb{R})}.
\end{equation*}%
Here $\eta $ depends only on $d$, $m_{d}$, and $(p_{1},\ldots ,p_{d})$, and $%
C$ depends only on $\beta $, $(m_{1},\ldots ,m_{d})$, and $(p_{1},\ldots
,p_{d})$. The estimate is uniform in $\lambda $, $j$, and $\tau $.
\end{prop}

\begin{proof}
We first choose $N+1$ equally spaced grid points $b_0,\ldots,b_N$ in $%
[1/2,1] $, where 
\begin{equation*}
N=\bigl\lceil\lambda^{c_0/2}\bigr\rceil\qquad\text{and}\qquad b_\ell=\frac12+%
\frac{\ell}{2N},\quad 0\leq\ell\leq N.
\end{equation*}
For any $b\in[1/2,1]$, choose a nearest grid point $b_\ell$. Then, for
almost every $x\in Q_k$, 
\begin{equation*}
|S_b(\vec g)(x)-S_{b_\ell}(\vec g)(x)|\leq
2|b-b_\ell|\prod_{i=1}^{d}\|g_i\|_{L^\infty(Q_k^*)}\leq
\frac1N\prod_{i=1}^{d}\|g_i\|_{L^\infty(Q_k^*)}.
\end{equation*}
Recall that $Q_k=[k,k+1]$ and $Q_k^*=[k-1,k+2]$. The triangle inequality
therefore gives 
\begin{equation*}
\big\|\sup_{1/2\leq b\leq 1}|S_b(\vec g)|\big\|_{L^1(Q_k)}\leq\sum_{%
\ell=0}^{N}\|S_{b_\ell}(\vec g)\|_{L^1(Q_k)}+\frac{1}{N}\prod_{i=1}^{d}\|g_i%
\|_{L^\infty(Q_k^*)}\lesssim
\lambda^{-c_0/2}\prod_{i=1}^{d}\|g_i\|_{L^\infty(Q_k^*)}
\end{equation*}
where the last inequality follows by applying Proposition~\ref{prop PET}.
Fix a measurable function $b=b(x)$ taking values in $[1/2,1]$, we have 
\begin{equation*}
\big\|S_{b(\cdot)}(\vec g)\big\|_{L^\infty(Q_k)}\lesssim
\|g_1\|_{L^1(Q_k^*)}\prod_{i=2}^{d}\|g_i\|_{L^\infty(Q_k^*)}.
\end{equation*}
With $g_2,\ldots,g_d$ fixed, we apply the interpolation theorem to this
estimate and the preceding $L^1(Q)$ estimate to obtain 
\begin{equation*}
\big\|S_{b(\cdot)}(\vec g)\big\|_{L^2(Q_k)}\lesssim\lambda^{-c_0/4}\|g_1%
\|_{L^2(Q_k^*)}\prod_{i=2}^d\|g_i\|_{L^\infty(Q_k^*)}.
\end{equation*}
Summing over $k\in\mathbb{Z}$, we obtain 
\begin{align}
\big\|S_{b(\cdot)}(\vec g)\big\|_{L^2(\mathbb{R})}\lesssim\lambda^{-c_0/4}%
\|g_1\|_{L^2(\mathbb{R})}\prod_{i=2}^d\|g_i\|_{L^\infty(\mathbb{R})}.
\label{eq1}
\end{align}
Choose a nonnegative $\rho\in C_c^\infty((1/4,2))$ such that $\rho(u)\geq
1/u $ for $u\in[1/2,1]$. Then, for every $b\in[1/2,1]$, we have the
pointwise bound 
\begin{equation*}
|S_b(\vec g)(x)|\leq\int_{\mathbb{R}}\rho(u)\prod_{i=1}^d
|g_i(x-P_i(u))|\,du.
\end{equation*}
Theorem \ref{thm GLZ} now gives 
\begin{align}
\|S_{b(\cdot)}(\vec g)\|_{L^s(\mathbb{R})}\lesssim\prod_{i=1}^d
\|g_i\|_{L^{a_i}(\mathbb{R})},  \label{eq2}
\end{align}
whenever $1<a_i<\infty$, $1/s=\sum_{i=1}^d 1/a_i$, and $s>1/2$. The implicit
constant is independent of the measurable choice $b(\cdot)$. Choose $%
\theta\in(0,1)$ sufficiently small and define $a_i$ and $s$ by 
\begin{equation*}
\frac{1}{a_i}=\frac{1/p_i-(\theta/2)1_{\{i=1\}}}{1-\theta},\qquad\frac{1}{s}=%
\frac{1/r-\theta/2}{1-\theta}.
\end{equation*}
Applying the multilinear complex interpolation theorem \cite[Theorem~3.2]{GM}
to \eqref{eq1} and \eqref{eq2}, we obtain 
\begin{equation*}
\|S_{b(\cdot)}(\vec g)\|_{L^r(\mathbb{R})}\lesssim\lambda^{-\theta
c_0/4}\prod_{i=1}^d \|g_i\|_{L^{p_i}(\mathbb{R})}.
\end{equation*}
Let $(b_\ell)_{\ell\geq1}$ be a dense sequence in $[1/2,1]$. For each $%
M\geq1 $, choose a measurable $b=b(x)$ from $\{b_1,\ldots,b_M\}$ so that 
\begin{equation*}
|S_{b(\cdot)}(\vec g)(x)|=\max_{1\leq\ell\leq M}|S_{b_\ell}(\vec g)(x)|.
\end{equation*}
The preceding estimate is uniform in this choice. Since $b\mapsto S_b(\vec
g)(x)$ is continuous for almost every $x$, the monotone convergence theorem
gives the desired supremum estimate as $M\to\infty$.
\end{proof}

In the next section, we prove Theorem~\ref{main thm} by combining a dyadic
annular decomposition of $T_{\beta,\gamma_d,\varepsilon}$ with the decay
estimate in Proposition~\ref{prop Sb}.

\section{proof of theorem \protect\ref{main thm} and corollary \protect\ref%
{cor:polynomial-curves}}

We now prove Theorem \ref{main thm}. We first treat the case in which all $%
m_{i}$ are positive integers.

The definition of the $1$-variation gives 
\begin{equation*}
V_{1}\bigl(S_{b}(\vec{g}):b\in \lbrack 1/2,1]\bigr)(x)\leq
\int_{1/2}^{1}\prod_{i=1}^{d}|g_{i}(x-P_{i}(u))|\,\frac{du}{u}.
\end{equation*}%
Applying theorem \ref{thm GLZ} once more, we obtain 
\begin{equation*}
\bigl\|V_{1}\bigl(S_{b}(\vec{g}):b\in \lbrack 1/2,1]\bigr)\bigr\|_{L^{r}(%
\mathbb{R})}\lesssim \prod_{i=1}^{d}\Vert g_{i}\Vert _{L^{p_{i}}(\mathbb{R}%
)}.
\end{equation*}%
For any finite partition $1/2\leq b_{0}<\cdots <b_{N}\leq 1$, we have 
\begin{equation*}
\sum_{l=1}^{N}\big|S_{b_{l}}(\vec{g})-S_{b_{l-1}}(\vec{g})\big|^{q}\leq \big(%
2\sup_{b\in \lbrack 1/2,1]}|S_{b}(\vec{g})|\big)^{q-1}\sum_{l=1}^{N}\big|%
S_{b_{l}}(\vec{g})-S_{b_{l-1}}(\vec{g})\big|.
\end{equation*}%
Thus, for $q>1$ 
\begin{equation*}
V_{q}\big(S_{b}(\vec{g}):\frac{1}{2}\leq b\leq 1\big)\leq \big(2\sup_{b\in
\lbrack 1/2,1]}|S_{b}(\vec{g})|\big)^{1-\frac{1}{q}}V_{1}\big(S_{b}(\vec{g}):%
\frac{1}{2}\leq b\leq 1\big)^{\frac{1}{q}}.
\end{equation*}%
Taking the $L^{r}$ quasi-norm and applying H\"{o}lder's inequality and
Proposition \ref{prop Sb}, we obtain 
\begin{equation*}
\begin{aligned} \bigl\|V_q\bigl(S_b(\vec
g):b\in[1/2,1]\bigr)\bigr\|_{L^r(\mathbb R)}&\lesssim\bigl\|\sup_{1/2\leq
b\leq1}|S_b(\vec g)|\bigr\|_{L^r(\mathbb
R)}^{1-1/q}\,\bigl\|V_1\bigl(S_b(\vec
g):b\in[1/2,1]\bigr)\bigr\|_{L^r(\mathbb R)}^{1/q}\\
&\lesssim\lambda^{-\eta(1-1/q)}\prod_{i=1}^{d}\|g_i\|_{L^{p_i}(\mathbb R)}.
\end{aligned}
\end{equation*}%
Next, for each $j\geq 0$ and $\sigma \in \{-1,1\}$, define 
\begin{equation*}
F_{j,\sigma }(b,x)=\int_{2^{-j-1}}^{2^{-j}b}e^{it^{-\beta
}}\prod_{i=1}^{d}f_{i}\bigl(x-\gamma _{d,i}(\sigma t)\bigr)\,\frac{dt}{t}%
,\qquad \frac{1}{2}\leq b\leq 1.
\end{equation*}%
The normalized polynomials have $\tau _{i}=1$ when $\sigma =1$ and $\tau
_{i}=\sigma _{i}$ when $\sigma =-1$. Set $g_{i}(x)=f_{i}(2^{-m_{1}j}x)$ and $%
\lambda =2^{\beta j}$. The change of variables $t=2^{-j}u$ gives 
\begin{equation*}
F_{j,\sigma }(b,2^{-m_{1}j}x)=S_{b}(\vec{g})(x).
\end{equation*}%
Combining the preceding estimates, we obtain 
\begin{align}
\bigl\|V_{q}\bigl(F_{j,\sigma }(b,\cdot ):b\in \lbrack 1/2,1]\bigr)\bigr\|%
_{L^{r}}& =2^{-m_{1}j/r}\bigl\|V_{q}\bigl(S_{b}(\vec{g}):b\in \lbrack 1/2,1]%
\bigr)\bigr\|_{L^{r}}  \notag \\
& \lesssim 2^{-\beta \eta (1-1/q)j}\prod_{i=1}^{d}\Vert f_{i}\Vert
_{L^{p_{i}}}.  \label{eq3}
\end{align}%
Finally, set $b_{j}(\varepsilon ):=\min \{1,\max \{1/2,2^{j}\varepsilon \}\}$%
. Then 
\begin{equation*}
\begin{aligned} T_{\beta,\gamma_d,\varepsilon}(\vec
f)(x)&=\sum_{\sigma=\pm1}\int_\varepsilon^1e^{it^{-\beta}}\prod_{i=1}^d
f_i\bigl(x-\gamma_{d,i}(\sigma t)\bigr)\,\frac{dt}{t}\\
&=\sum_{\sigma=\pm1}\sum_{j\geq
0}\big[F_{j,\sigma}(1,x)-F_{j,\sigma}(b_j(\varepsilon),x)\big]. \end{aligned}
\end{equation*}%
Since $b_{j}(\varepsilon )$ is nondecreasing in $\varepsilon $, the triangle
inequality for $V_{q}$ gives 
\begin{equation*}
V_{q}\bigl(T_{\beta ,\gamma _{d},\varepsilon }(\vec{f}):0<\varepsilon <1%
\bigr)\leq \sum_{\sigma =\pm 1}\sum_{j\geq 0}V_{q}\bigl(F_{j,\sigma
}(b,\cdot ):b\in \lbrack 1/2,1]\bigr).
\end{equation*}%
Setting $s_{0}=\min \{1,r\}$ and applying the estimate (\ref{eq3}), we
finally obtain 
\begin{equation*}
\begin{aligned} \bigl\|V_q\bigl(T_{\beta,\gamma_d,\varepsilon}(\vec
f):0<\varepsilon<1\bigr)\bigr\|_{L^r}^{s_0}
&\leq\sum_{\sigma=\pm1}\sum_{j\geq0}\bigl\|V_q\bigl(F_{j,\sigma}(b,\cdot):b%
\in[1/2,1]\bigr)\bigr\|_{L^r}^{s_0}\\
&\lesssim\sum_{j\geq0}2^{-\beta\eta(1-1/q)js_0} \prod_{i=1}^d
\|f_i\|_{L^{p_i}}^{s_0}\lesssim\prod_{i=1}^d \|f_i\|_{L^{p_i}}^{s_0}.
\end{aligned}
\end{equation*}%
Taking the $s_{0}$-th root completes the proof in the integer case.

For any rational exponents $0<m_{1}<\cdots <m_{d}$, let $L$ be a common
denominator of $m_{1},\ldots ,m_{d}$. Let $n_{i}=Lm_{i}$, then $n_{1},\ldots
,n_{d}$ are positive integers. Making the changes of variables, we have 
\begin{equation*}
\begin{aligned} T_{\beta,\gamma_d,\varepsilon}(\vec
f)=L\int_{\varepsilon^{1/L}}^{1}e^{i
s^{-L\beta}}\prod_{i=1}^{d}f_i(x-s^{n_i})\frac{ds}{s}+L\int_{%
\varepsilon^{1/L}}^{1}e^{i
s^{-L\beta}}\prod_{i=1}^{d}f_i(x-\sigma_is^{n_i})\frac{ds}{s} \end{aligned}
\end{equation*}%
where $\sigma _{i}\in \{-1,1\}$ is a fixed sign. Applying the estimate for
integer exponents to each of the two terms above, we obtain the theorem \ref%
{main thm}.

It remains to prove the corollary \ref{cor:polynomial-curves}. Without loss
of generality, we may assume that $P_{i}(0)=0$ for every $i$. If all the $%
P_{i}$ vanish identically, the result follows by integration by parts on
each dyadic annulus and summation over $j$. We may therefore assume that at
least one $P_{i}$ is nonzero. Let $m$ be the smallest order of vanishing at
the origin among the nonzero $P_{i}$. Relabel the inputs so that 
\begin{equation*}
P_{1}(t)=c_{1}t^{m}+O(t^{m+1})\qquad c_{1}\neq 0,\,|t|\leq 1.
\end{equation*}%
For $j\geq 0$ and $\sigma \in \{-1,1\}$, set 
\begin{equation*}
\widetilde{P}_{i,j,\sigma }(u)=\frac{P_{i}(\sigma 2^{-j}u)}{A2^{-mj}},
\end{equation*}%
where $A$ is chosen sufficiently large that $|\widetilde{P}_{i,j,\sigma
}(u)|\leq 1$ for $u\in \lbrack 1/2,1]$. There exists an integer $j_{0}\geq 0$
such that, for every $j\geq j_{0}$ and $\sigma \in \{-1,1\}$, 
\begin{equation*}
|\widetilde{P}_{1,j,\sigma }^{\prime }(u)|\geq \frac{m|c_{1}|}{A2^{m}}%
>0,\qquad u\in \lbrack 1/2,1].
\end{equation*}%
Moreover, $\widetilde{P}_{1,j,\sigma }$ is strictly monotone on this
interval. Define $S_{b}$ as in Proposition~\ref{prop Sb}, with $P_{i}$
replaced by $\widetilde{P}_{i,j,\sigma }$. For $x\in Q_{k}$, all the shifted
arguments lie in $Q_{k}^{\ast }$. The change of variables $y=x-\widetilde{P}%
_{1,j,\sigma }(u)$ therefore gives 
\begin{equation*}
\begin{aligned} |S_{b(x)}(\vec
g)(x)|&\leq2\prod_{i=2}^{d}\|g_i\|_{L^\infty(Q_k^*)}\int_{1/2}^{1}|g_1(x-%
\widetilde P_{1,j,\sigma}(u))|\,du\\
&\lesssim\|g_1\|_{L^1(Q_k^*)}\prod_{i=2}^{d}\|g_i\|_{L^\infty(Q_k^*)}.
\end{aligned}
\end{equation*}%
The proof of proposition~\ref{prop PET} still applies to these polynomial
shifts. Theorem~\ref{thm GLZ} applies only to monomial shifts. We therefore
replace it with the following estimate. For $s>1$ and $1/s=%
\sum_{i=1}^{d}1/a_{i}$, Minkowski's and H\"{o}lder's inequalities give 
\begin{equation*}
\begin{aligned} \Big\|\int_{1/2}^{1}\prod_{i=1}^d|g_i(\,\cdot-\widetilde
P_{i,j,\sigma}(u))|\,\frac{du}{u}\Big\|_{L^s} &\leq
\int_{1/2}^{1}\Big\|\prod_{i=1}^d |g_i(\,\cdot-\widetilde
P_{i,j,\sigma}(u))|\Big\|_{L^s}\frac{du}{u}\\ &\leq (\log
2)\prod_{i=1}^d\|g_i\|_{L^{a_i}}. \end{aligned}
\end{equation*}%
For the finitely many annuli with $0\leq j<j_{0}$, we use $V_{q}\leq V_{1}$
and apply the preceding estimate with $s=r$ and $a_{i}=p_{i}$.

For $j\geq j_0$, set $g_i(x)=f_i(A2^{-mj}x)$ and $\lambda=2^{\beta j}$.
Since $r>1$, we may choose the interpolation parameter sufficiently small
that $s>1$. The rest of the proof proceeds exactly as in the proof of
Theorem~\ref{main thm}.

\end{document}